\documentstyle [12pt]{article}

\newfam\msbfam
\font\twelmsb=msbm10 at 12pt
\font\tenmsb=msbm10
\font\sevenmsb=msbm10 at 7pt
\font\fivemsb=msbm10 at 5pt

\textfont\msbfam=\twelmsb
\scriptfont\msbfam=\sevenmsb
\scriptscriptfont\msbfam=\fivemsb

\def\Bbb{\fam\msbfam\tenmsb}

\def\R{{\Bbb R}}
\def\Z{{\Bbb Z}}

\def\KK{{\cal K}}

\def\OO{{\cal O}}

\def\WW{{\cal W}}

\def\qed {\hfill\vrule height6pt width6pt depth0pt \bigskip}

\def\map{\longrightarrow}
\def\textmap#1{\mathop{\vbox{\ialign{
                                ##\crcr
    ${\scriptstyle\hfil\;\;#1\;\;\hfil}$\crcr
    \noalign{\kern-1pt\nointerlineskip}
    \rightarrowfill\crcr}}\;}}

\def\textlmap#1{\mathop{\vbox{\ialign{
                                ##\crcr
    ${\scriptstyle\hfil\;\;#1\;\;\hfil}$\crcr
    \noalign{\kern-1pt\nointerlineskip}
    \leftarrowfill\crcr}}\;}}

\newfam\meuffam
\font\tenmeuf=eufm10
\font\sevenmeuf=eufm7
\font\fivemeuf=eufm5

\textfont\meuffam=\tenmeuf
\scriptfont\meuffam=\sevenmeuf
\scriptscriptfont\meuffam=\fivemeuf

\def\germ{\fam\meuffam\tenmeuf}

\def\picture#1by#2(#3){
\vbox to #2 {
  \hrule width #1 height 0pt depth 0pt \vfill \special{picture #3}}
}

\def\scaledpicture#1by#2(#3scaled#4){{
\dimen0=#1  \dimen1=#2
\divide\dimen0 by 1000 \multiply\dimen0 by #4
\divide\dimen1 by 1000 \multiply\dimen1 by #4
\picture \dimen0 by \dimen1 (#3 scaled #4)}}
\def\dfigure#1by#2(#3scaled#4offset#5:#6)
  {\medskip
   \vglue 2mm minus 2mm
   $$
     \hbox{
       \hglue#5
       {\scaledpicture #1 by #2 (#3 scaled #4)}
     }
   $$
   \par\goodbreak
   \vglue 2mm minus 2mm
   \medskip}

\begin{document}
\def\Pr{{\rm Pr}}
\def\tr{{\rm Tr}}
\def\End{{\rm End}}
\def\Pic{{\rm Pic}}
\def\NS{{\rm NS}}
\def\deg{{\rm deg}}
\def\Hom{{\rm Hom}}
\def\Herm{{\rm Herm}}
\def\Vol{{\rm Vol}}
\def\pf{{\bf Proof: }}
\def\id{{\rm id}}
\def\i{{\germ i}}
\def\im{{\rm im}}
\def\rk{{\rm rk}}
\def\h{{\bf H}}
\def\dv{\bar\partial}
\def\dva{\bar\partial_A}
\def\da{\partial_A}
\def\p{\partial\bar\partial}
\def\pa{\partial_A\bar\partial_A}
\def\Dr{\hskip 4pt{\not}{D}}
\newtheorem{sz}{Satz}
\newtheorem{th}[sz]{Theorem}
\newtheorem{pr}[sz]{Proposition}
\newtheorem{re}[sz]{Remark}
\newtheorem{co}[sz]{Corollary}
\newtheorem{dt}[sz]{Definition}
\newtheorem{lm}[sz]{Lemma}
\title{Seiberg-Witten Invariants and the Van De Ven Conjecture}
\author{Christian Okonek$^*$\and Andrei Teleman\thanks{Partially supported
by: AGE-Algebraic Geometry in Europe, contract No ERBCHRXCT940557 (BBW
93.0187), and
by  SNF, nr. 21-36111.92} }
\date{February, 8-th 1995}
\maketitle

The purpose of this note is to give a short, selfcontained proof of the
following
result:
\begin{th}
A  complex surface which is  diffeomeorphic to a rational surface is rational.
\end{th}
This result  has been announced by R. Friedman and Z. Qin [FQ].
Whereas their proof uses Donaldson theory and vector bundles techniques, our
proof uses the new Seiberg-Witten invariants [W], and the interpretation of
these invariants in terms of stable pairs [OT].

Combining the theorem above with the results of [FM], one obtains a proof of
the
Van de Ven conjecture [V]:
\begin{co}
The Kodaira dimension of a complex surface is a differential invariant.
\end{co}

\pf (of the Theorem) It suffices to prove the theorem for algebraic
surfaces [BPV].
Let $X$ be an algebraic surface of non-negative Kodaira dimension, with
$\pi_1(X)=\{1\}$ and $p_g(X)=0$. We may suppose that $X$ is the blow up in $k$
{\sl distinct} points of its minimal model $X_{\min}$. Denote the
contraction to the
minimal model by
$\sigma:X\map X_{\min}$, and the exceptional divisor by $E=\sum\limits_{i=1}^k
E_i$.

Fix an ample divisor $H_{\min}$ on $X_{\min}$, a sufficiently large integer
$n$, and
let $H_n:=\sigma^*(n H_{\min})-E$ be the associated polarization of $X$.

For every subset $I\subset\{1,\dots,k\}$ we put
$E_I:=\sum\limits_{i\in I} E_i$, and
$L_I:=2[E_I]-[K_X]$, where $K_X$ is
a canonical divisor. Clearly $L_I=[E_I]-[E_{\bar I}]-\sigma^*([K_{\min}])$,
where
$\bar I$ denotes the complement of $I$ in $\{1,\dots,k\}$. The cohomology
classes
$L_I$ are almost canonical classes in the sense of [OT].  Now choose
a K\"ahler metric $g_n$ on $X$ with K\"ahler class
$[\omega_{g_n}]=c_1(\OO_X(H_n))$. Since $[\omega_{g_n}]\cdot L_I<0$ for
sufficiently large
$n$, the main result of [OT] identifies the Seiberg-Witten moduli space
$\WW_X^{g_n}(L_I)$ with the union of all complete linear systems $|D|$
corresponding to effective divisors $D$ on $X$ with
$c_1(\OO_X(2D-K_X))=L_I$.

Since $H^2(X,\Z)$ has no 2-torsion, and $q(X)=0$, there is only one such
divisor, $D=E_I$. Furthermore, from
$h^1(\OO_X(E_I)|_{E_I})=0$, and the smoothness criterion in [OT], we obtain:
$$\WW_X^{g_n}(L_I)=\{E_I\}, $$
i.e. $\WW_X^{g_n}(L_I)$ consists of a single smooth point.
The corresponding Seiberg-Witten invariants are therefore odd:
$n_{L_I}^{g_n}=\pm 1$.

Consider now the positive cone  $\KK:=\{\eta\in H^2_{\rm DR}(X)|\ \eta^2>0\}$;
using the Hodge index theorem, the fact that $K_{\min}$ is cohomologically
nontrivial, and $K_{\min}^2\geq 0$, we see that
$\KK$ splits as a disjoint union of two components  $\KK_{\pm}:=\{\eta\in\KK|\
\pm\eta\cdot\sigma^*(K_{\min})>0\}$. Clearly $[\omega_{g_n}]$ belongs to
$\KK_+$.

Let $g$ be an {\sl arbitrary} Riemannian metric on $X$, and let $\omega_g$ be
a $g$-selfdual closed 2-form on $X$ such that $[\omega_g]\in\KK_+$.

For a fixed $I\subset\{1,\dots,k\}$, we denote by $L_I^{\bot}\subset\KK_+$ the
wall associated with $L_I$, i.e. the subset of classes $\eta$ with $\eta\cdot
L_I=0$.\\ \\
{\bf Claim:} The rays $\R_{>0}[\omega_g]$,  $\R_{>0}[\omega_{g_n}]$ belong
either to the same component of $\KK_+\setminus L_I^{\bot}$ or to the same
component of   $\KK_+\setminus L_{\bar I}^{\bot}$. \\

Indeed, since $[\omega_{g_n}]\cdot L_I<0$ and $[\omega_{g_n}]\cdot L_{\bar
I}<0$,
we just have to exclude that
$$[\omega_{g}]\cdot L_I\geq 0 \ \ \ \ {\rm and}\ \ \ \ [\omega_{g}]\cdot
L_{\bar
I}\geq 0 .\eqno{(*)}$$

Write $[\omega_g]=\sum\limits_{i=1}^{k}\lambda_i[E_i]+\sigma^*[\omega]$,
for some
class
$[\omega]\in H^2_{\rm DR}(X_{\min})$; then
$[\omega]^2>\sum\limits_{i=1}^k\lambda_i^2$, and $[\omega]\cdot
K_{\min}>0$, since
$\omega_g$ was chosen such that its cohomology class belongs to $\KK_+$. The
inequalities $(*)$ can now be written as
$$-\sum\limits_{i\in I}\lambda_i+\sum\limits_{j\in\bar I}\lambda_j-
[\omega]\cdot K_{\min}\geq 0 \ \ \ {\rm and}\ \ \
-\sum\limits_{j\in\bar I}\lambda_j+\sum\limits_{i\in I}\lambda_i-
[\omega]\cdot K_{\min}\geq 0,$$
and we obtain the contradiction $[\omega]\cdot K_{\min}\leq 0$. This proves the
claim.

\dfigure 80mm by 165mm (kegel scaled 700 offset 1mm:)

We know already that the mod 2 Seiberg-Witten invariants $n^{g_n}_{L_I}$(mod 2)
and  $n^{g_n}_{L_{\bar I}}$(mod 2) are nontrivial for the special metric $g_n$.
Since  the invariants $n^{g}_{L_{I}}$(mod 2) and $n^{g}_{L_{\bar I}}$(mod 2)
depend only on the chamber of the ray $\R_{>0}[\omega_g]$  with repect to the
wall $L_{I}^{\bot}$,  respectively $L_{\bar I}^{\bot}$ (see [W], [KM]), at
least
one of the invariants associated with the metric $g$ must  be non-zero, too.

But any rational surface admits a Hodge metric with positive total scalar
curvature [H], and with respect to such a metric {\sl all} Seiberg-Witten
invariants
are trivial [OT].
\qed
\vspace{0.5cm}\\
\parindent0cm
\centerline {\Large {\bf Bibliography}}

\vspace{0.5cm}
[BPV] Barth, W., Peters, C., Van de Ven, A.: {\it Compact complex surfaces},
Springer Verlag (1984)

[FM] Friedman, R., Morgan, J.W.: {\it Smooth 4-manifolds and Complex Surfaces},
Springer Verlag  3. Folge, Band 27, (1994)

[FQ]  Friedman, R., Qin, Z.: {\it On complex surfaces diffeomorphic to
rational surfaces}, Preprint (1994)

[H] Hitchin, N.: {\it  On the curvature of rational surfaces}, Proc. of Symp.
in Pure Math., Stanford, Vol. 27 (1975)

[KM] Kronheimer, P., Mrowka, T.: {\it The genus of embedded surfaces in the
projective plane}, Preprint (1994)

[OT] Okonek, Ch.; Teleman A.: {\it The Coupled Seiberg-Witten Equations,
Vortices, and Moduli Spaces of Stable Pairs}, Preprint, January, 13-th 1995

[W] Witten, E.: {\it Monopoles and four-manifolds}, Mathematical Research
Letters 1, 769-796 (1994)

[V] Van de Ven, A,: {\it On the differentiable structure of certain algebraic
surfaces}, S\'em. Bourbaki ${\rm n}^o$ 667, Juin (1986)
\vspace{1cm}\\
Authors addresses:\\
\\
Mathematisches Institut, Universit\"at Z\"urich,\\
Winterthurerstrasse 190, CH-8057 Z\"urich\\
e-mail:okonek@math.unizh.ch

\ \ \ \ \ \  \ \ \ teleman@math.unizh.ch

\end{document}